\documentclass{amsart}
\usepackage{latexsym}
\usepackage{amssymb}
\usepackage{amsmath}

\newtheorem{thm}{Theorem}[section]
\newtheorem{prop}[thm]{Proposition}
\newtheorem{lem}[thm]{Lemma}
\newtheorem{cor}[thm]{Corollary}
\def\dem{\par\noindent \emph{Proof.} \ \,}
\def\edem{\rm\hfill$\Box$\bigskip}
\def\vhi{\varphi}
\def\eps{\varepsilon}

\def\bs{\backslash}
\def\aut{\operatorname{Aut}}

\def\m{^{-1}}

\def\inn{\operatorname{Inn}}

\def\ml{\mathcal L}
\def\mr{\mathcal R}

\newcommand{\cref}[1]{Corollary~\ref{#1}}
\newcommand{\lref}[1]{Lemma~\ref{#1}}
\newcommand{\pref}[1]{Proposition~\ref{#1}}

\newcommand{\secref}[1]{Section~\ref{#1}}

\begin{document}
\title[Buchsteiner loops]
{Buchsteiner loops: associators and constructions}
\author {Ale\v s Dr\'apal}
\author{Michael Kinyon}

\address{Dept.~of Mathematics \\ Charles University \\ Sokolovsk\'a 83 \\ 186
75 Praha 8, Czech Rep.}
\email{drapal@karlin.mff.cuni.cz}
\address{Dept.~of Mathematics \\ University of Denver \\2360 S. Gaylord St. \\
Denver, Colorado 80208, U.S.A.}
\email{mkinyon@math.du.edu}

\thanks{ The first
author (A. Dr\'apal) was supported by institutional grant MSM 0021620839 and
by Grant Agency of Charles University, Grant 444/2004. This paper was written
while he was a Fulbright Research Scholar at University of
Wisconsin-Madison.
}

\begin{abstract}
Let $Q$ be a Buchsteiner loop. We describe the associator calculus in
three variables, and show that $|Q| \ge 32$ if $Q$ is not conjugacy closed.
We also show that $|Q| \ge 64$ if there exists $x \in Q$ such that $x^2$ is
not in the nucleus of $Q$. Furthermore, we describe a general construction
that yields all proper Buchsteiner loops of order $32$. Finally, we produce
a Buchsteiner loop of order $128$ that is nilpotency class $3$ and possesses
an abelian inner mapping group.
\end{abstract}

\maketitle

Buchsteiner loops are those loops that satisfy the Buchsteiner law
$$x\bs (xy\cdot z) = (y \cdot zx)/x.$$
Their study was initiated by Hans Hoenig Buchsteiner \cite{hhb}. His paper
left many problems open, some of which were recently solved \cite{buch}.
In particular we know now that the nucleus $N = N(Q)$ is a normal
subloop of every Buchsteiner loop $Q$ and that $Q/N$ is an abelian
group of exponent four. 

Buchsteiner loops are closely connected to conjugacy closed loops (CC loops).
A CC loop is conjugacy closed if and only if $Q/N$ is a boolean group
(i.e. a group of exponent two), by \cite{nuc}. Not every Buchsteiner loop
with $Q/N$ boolean needs to be conjugacy closed (there are plenty of examples
now. Some of them appear in this paper, and many other can be
derived from the ring construction of \cite{bcc}.)

In every Buchsteiner loop $Q$ the mappings $L_{xy}\m L_x L_y$ and 
$R_{yx}\m R_x R_y$ are automorphisms of $Q$, by \cite{buch}, and this fact
effects the behaviour of the associators $[x,y,z] = (x\cdot yz)\bs (xy \cdot z)$.
The group $Q/A(Q)$ acts upon $N = N(Q)$ (that always holds when $Q/N$ is a group
since then $A(Q) \le Z(N(Q))$, by \cite{dia}. Here $A(Q)$ denotes the
least normal subloop $A \unlhd Q$ such that $Q/A$ is a group. If $Q/N$ is 
a group, then $A(Q)$ coincides with the subgroup generated by all associators
$[x,y,z]$, by \cite{pac}). By translating the automorphism behaviour
of $L_{xy}\m L_x L_y$ into relations between associators we get that
$[x,y,uv] = [x,y,u]^v[x,y,v]$ for all $x,y,u,v \in Q$. Note that $n^v$ is
defined as $v\bs (nv)$, for all $v \in Q$ and $n \in N$. 

If $Q$ is a loop such that $Q/N$ is a group, then one can code the Buchsteiner 
identity as 
$$ [x,y,z]^x = [y,z,x]\m \text{\, for all \,} x,y,z \in Q.$$

The cyclic shift expressed by this action implies that in every Buchsteiner
loop we have
$$[x,y,uv] = [x,y,u]^v[x,y,v], [x,uv,y]=[x,u,y]^v[x,v,y],
[uv,x,y] = [u,x,y]^v [v,x,y]$$
for all $x,y,u,v \in Q$. If $Q/N$ is a group, then an associator $[x,y,z]$
depends only upon the ordered triple $(xN,yN,zN)$, by \cite{dia}. If $Q$
is a Buchsteiner loop, then $Q/N$ is an abelian group, and so 
we have $[x,y,uv] = [x,y,vu] = [x,y,v]^u[x,y,u]
= [x,y,u][x,y,v]^u$. Similar relations clearly hold for the other two positions
too. All the
facts above are exposed in \cite{buch} in detail, and we shall use them freely
in this paper.

In \secref{1} we shall develop the associator calculus in three variations,
building upon the identities that were established in \cite{buch}.
Sections~\ref{2}--\ref{4} are mainly concerned with the proof
that proper Buchsteiner loops are of order at least $32$ (by a \emph{proper}
Buchsteiner loop we understand a Buchsteiner loop that is not
conjugacy closed). Buchsteiner loops $Q$ such that
$Q/N$ is not boolean are necessarily proper, and for them we show
that $|Q| \ge 64$. In \secref{4} we will observe that nilpotent proper
Buchsteiner loops have to be of nilpotency class at least $3$, a result
that appears also in \cite{bcc}.

In \secref{6} we will show that a proper Buchsteiner loop of order 32 
really exists, and 
that all such loops can be obtained by a general construction that
doubles the size of a loop. This construction is described in \secref{5}.
The starting loop must be a Buchsteiner loop, but not necessarilly 
a proper one.

In a loop $Q$ it is usual to denote by $L_x$ the left translation $y \mapsto xy$,
and by $R_x$ the right translation $y \mapsto yx$. The permutation group
generated by all $L_x$ and $R_x$ is known as the \emph{multiplication group},
and the stabilizer of the unit is called the \emph{inner mapping group};
we denote it by $\inn Q$.
It is well known that $\inn Q$ is generated by all mappings
$L(x,y) = L_{xy}\m L_xL_y$, $R(x,y) = R_{yx}\m R_x R_y$ and $T_x = R_x\m L_x$.
If $Q$ is of nilpotency class two, then the inner mapping group is abelian,
a result that goes back to Bruck \cite{con}. The converse is not true,
but the examples are not easy to find. Up to now there has been
published only one example, by Cs\"org\H o \cite{csc}. Her example has
128 elements, was constructed indirectly by means of group transversals
(see also \cite{csa}),
and does not belong to any of known specific loop classes. 
In Sections~\ref{7} and~\ref{8} we construct a proper Buchsteiner loop $Q$
of order $128$ with $\inn Q$ abelian. This loop is different
from the construction of \cite{csc}, and is necessarily
of nilpotency class three since the Buchsteiner loops of nilpotency
class two are conjugacy closed. \secref{7} is concerned with general properties
of Buchsteiner loops that have abelian inner mapping groups,
and \secref{8} contains the construction. Note that (left) conjugacy 
closed loops with abelian inner mapping groups are always of nilpotency
class at most two, by \cite{csd}.

\section{Associator identities}\label{1}

Let $Q$ be a Buchsteiner loop. Then $[z\m,x,y]^{z\m} = [x,y,z\m]\m$,
which we write as $[x,y,z\m]^z = [z\m,x,y]\m$. Therefore $1 = [x,y,z\m z]
= [x,y,z\m]^z [x,y,z] = [z\m,x,y]\m [x,y,z]$. Hence $[x,y,z] = [z\m,x,y]$.
This identity comes from \cite{buch}, where many further similar calculations
have been performed. The next two lemmas give a selection of them. 

\begin{lem}\label{11}
Let $Q$ be a Buchsteiner loop with elements $x$, $y$ and $z$. Put
$s = [x^2,y,z]$. Then $s = s^x = s^y = s^z$, $s^2 =1$, and $s = [y,z,x^2] = [z,x^2,y]$.
Furthermore, each of $x^2$, $y^2$ and $z^2$ centralizes $[x,y,z]$ (e.g.
$[x,y,z]^{x^2} = [x,y,z]$ etc.).
\end{lem}

\begin{lem}\label{12}
Let $Q$ be a Buchsteiner loop with elements $x$, $y$ and $z$. Put 
$u = [x,y\m,z]$. Then:
\begin{enumerate}
\item[(i)] $[x,y,z]^2 = [y,z,x]^2 = [z,x,y]^2 = u^2$;
\item[(ii)] $[x,y,z][y,z,x] = [z,x,y][x,y\m,z]$; 
\item[(iii)] $u = [x,y\m,z]= [y,z\m,x] = [z,x\m,y]$; and
\item[(iv)] $u^x = [z,x,y]\m$, $u^y = [x,y,z]\m$, and $u^z = [z,x,y]\m$.
\end{enumerate}
\end{lem}

The fact that $Q/N$ is of exponent four, in every Buchsteiner loop $Q$, means
that $[x^4,y,z] = 1$ and $[x\m,y,z] = [x^3,y,z]=[x^2,y,z][x,y,z]$, 
for all $x,y, z\in Q$.

\begin{prop}\label{13}
Let $Q$ be a Buchsteiner loop with elements $x$, $y$ and $z$. Put
$u = [x,y\m, z]$, $s_x = [x^2,y,z]$, $s_y = [y^2,z,x]$ and $s_z = [z^2,y,x]$.
Then $s_x = u^x u$, $s_y = u^yu$ and $s_z = u^z u$. Furthermore,
$$ s_x s_y s_z = 1,\ u^xu^yu^z = u^{-3},\  u^{xyz} = u\m, 
\text{ and \,} (u^2)^x = (u^2)^y = (u^2)^z = u^{-2}.$$
\end{prop}
\dem Each of $s_x$, $s_y$ and $s_z$ is of exponent two, and
$[z,x,y]^{x^2} = [z,x,y]$, 
by \lref{11}. From points (iii) and (iv) of \lref{12} we 
can compute $u^x u$ as $
[z,x\m,y]^x[z,x\m,y] = [z,x,y]\m[z,x,y][z,x^2,y] = s_x$, and the 
identities $u^y u = s_y$ and $u^z u = s_z$ can be proved similarly. 
Points (i) and (iv) of \lref{12} yield $(u^2)^x = (u^2)^y = (u^2)^z = u^{-2}$.
Point (ii) of the lemma can be written as $1 = u^y u^z (u^x)\m u$, and $(u^x)\m$ can
be replaced by $(u^x)(u^x)^{-2} = u^x u^2$. This means $u^xu^yu^z = u^{-3}$,
and so $s_xs_ys_z = u^xuu^yuu^zu = 1$. Finally, $u^{xyz} = u^{yxz} =
([x,y,z]\m)^{xz} = [y,z,x]^z = [y,z\m,x]\m =u\m$ (recall that $z^2$ centralizes
$[y,z,x]$, by \lref{11}).
\edem

Note that $[x,y,z]^{xy} = ([y,z,x]\m)^y = [z,x,y]$, for all $x,y,z \in Q$.
This gives $[x,x,y] = [x,x,y]^{x^2} = [y,x,x]$. We shall now prove some
further facts that involve only two variables. Most of the equalities can
be found in \cite{buch}, but we shall prove them here, in order to keep
the interface with \cite{buch} limited. (There are usually many ways
how one can obtain an identity. \pref{13} can be always used when an
associator is conjugated by a composition of its arguments. So we
can also get $[x,y,z]^{xy}$ as $(u\m)^{yxy} = (u\m)^x = [z,x,y]$.)

\begin{lem}\label{14}
Let $Q$ be a Buchsteiner loop with elements $x$ and $y$. Put
$u = [x,y,x]$ and $v = [x,x,y]$. Then
\begin{gather*}
u^y = u\m,\ v^y = v\m,\  u^x = v\m,\  v^x = u\m, 
\\ u^2 = v^2 = [y,x,y]^2 = [y,y,x]^2
\text{\, and \,} uv\m = vu\m = [x^2,x,y].
\end{gather*}
Furthermore, $[x,y^2,x] = 1$ and $[x,x,x]^y = [x,x,x](uv)\m$.
\end{lem}
\dem We have $[x,y^2,x] = [y^2,x,x]$, by \lref{11}, and $[y^2,x,x]$ is 
equal to
$[y,x,x]^y[y,x,x] = [x,x,y]\m[y,x,x] = 1$. Hence $u = [x,y\m,x] = [x,y,x]
[x,y^2,x]$, and the equalities $u^y = u\m$ and $u^x = v\m$
follow from \pref{13}. By \lref{11}, $u^2 = v^2$ (and so $uv\m = vu\m$),
and both of $u$ and $v$
are centralized by both of $x^2$ and $y^2$. Thus $u^x = v\m$ yields
$v^x = u\m$. Clearly, $[x^2,x,y] = v^xv = u\m v$, and $v[x,x,x]^y
= [x,x,xy] = [x,x,x]v^x = [x,x,x]u\m$. A similar argument can be used
to prove $[y,x,y]^2 = v^2$. Indeed,
$$ [y,x,y]v\m = [y,x,y]v^y =[yx,x,y] = v[y,x,y]^x = v[y,x,y]\m.$$
\edem

More results can be obtained by arguments similar to the one we used when
computing $[x,x,x]^y$:
\begin{lem}\label{15}
Let $Q$ be a Buchsteiner loop with elements $x$, $y$ and $z$. Then:
\begin{gather*}
[x,y,x]^{z^2} = [x,y,x],\ [x,x,y]^{z^2} = [x,x,y],\ 
[x^2,x,y]^z = [x^2,x,y] \text{\, and}\\
[x,y,x]^z = [x,y,x][x^2,y,z][z,x,y]^{-2} = [x,y,x][x^2,z,y][z,y,x]^{-2}.
\end{gather*}
\end{lem}
\dem Write $[xz^2,x,y]$ as $[x,x,y]^{z^2}[z^2,x,y]$
and as $[x,x,y][z^2,x,y]^x$. We know that $x$ centralizes
$[z^2,x,y]$, by \lref{11}, and hence $z^2$ centralizes $[x,x,y]$.
A similar argument shows that $z^2$ also centralizes $[x,y,x]$. 
Furthermore,
$[x^2,x,y][z,x,y] = [x^2,x,y][z,x,y]^{x^2} = [x^2z,x,y] =
[x^2,x,y]^z[z,x,y]$, and so $z$ centralizes $[x^2,x,y]$.

Write $[xz,x,y]$ as $[x,x,y]^z[z,x,y]$ and as $[x,x,y]
[z,x,y]^x = [x,x,y][z,x\m,y]\m$. Hence $[x,x,y]^z = [x,x,y]v$,
where $v = [z,x\m,y]\m[z,x,y]\m = [z,x^2,y][z,x,y]^{-2}$, and 
so $[x,y,x]^z = ([x,x,y][x^2,x,y])^z = [x,x,y][x^2,x,y]v 
= [x,y,x]v$, by the preceding parts of the proof and by \lref{14}.

Proceeding similarly, 
write $[xz,y,x]$ as $[x,y,x]^z[z,y,x]$ and as
$[x,y,x][z,y,x]^x = [x,y,x][x,z,y]\m$. Thus $[x,y,x]^z=
[x,y,x]w\m$, where $w = [z,y,x][x,z,y]
= [y,x,z][y,x\m,z] = [y,x^2,z][y,x,z]^2$, by points (ii) and (iii)
of \lref{12}.
\edem

To make complete our understanding of the associator calculus in
three variables we need to establish the relationship of associators
$[x,y,z]$ and $[y,x,z]$. 

\begin{prop}\label{16}
Let $Q$ be a Buchsteiner loop with elements $x$, $y$ and $z$. Then
\begin{align*} [x,y,z][z,y,x]\m &= [y,z,x][x,z,y]\m = [z,x,y][y,x,z]\m.
\end{align*}				
Denote this element by $a$, and put
$s_x = [x^2,y,z]$, $s_y = [y^2,z,x]$ and $s_z = [z^2,y,x]$. Then
\begin{align*}
a^2 = 1,\,a^x = a^y =a^z =a,\,[x,y,&z]^2 = [x,z,y]^2,\,
[x^2,y,z] = [x^2,z,y], \text{\, and}\\
[x,y,z][x,z,y]\m&=as_x= [y,z,x][z,y,x]\m\, \\ [y,z,x][y,x,z]\m &= as_y
= [z,x,y][x,z,y]\m
\text{\, and \,} \\ 
[z,x,y][z,y,x]\m&= as_z = [x,y,z][y,x,z]\m.
\end{align*}
\end{prop}
\dem 
We shall again use the equality $[x,y,z]^{xy} = [z,x,y]$. We obtain:
\begin{align*}
[xy,xy,z]=[x,xy,z]^y[y,xy,z] = [x,x,z]^y&[x,y,z]^{xy}[y,y,z]^x [y,x,z]=
\\ &[x,x,z]^y[y,y,z]^x[z,x,y][y,x,z], \text{ and} \\
[xy,xy,z] = [x,xy,z][y,xy,z]^x = [x,x,z]^y&[x,y,z][y,y,z]^x[y,x,z]^{yx} =
\\ &[x,x,z]^y [y,y,z]^x [x,y,z][z,y,x].
\end{align*}
By comparing the right hand sides we get $[x,y,z][y,x,z]\m = [z,x,y][z,y,x]\m$.
We also have $[x,y,z][y,x,z] = [z,x,y][z,y,x]$, since $[y,x,z]^2 = [z,y,x]^2$,
by point (i) of \lref{12}. Therefore,
$$[x,y,z][z,y,x]\m = [z,x,y][y,x,z]\m = [y,z,x][x,z,y]\m,$$
where the latter equality is an instance of the former one. Denote this 
element by $a$, as in the text of the proposition. Note that we have proved
that the leftmost term and the rightmost term coincide in all three
last equalities of the proposition.

Now, $[x^2,y,z] = [x,y,z][x,y,z]^2 = [x,y,z][y,z,x]\m$ equals 
$[z,y,x][x,z,y]\m = [z,y,x^2]= [x^2,z,y]$, and that immediately
yields $[x,y,z]^2 = [z,x,y]^2 = [z,y,x]^2 = [x,z,y]^2$, by \lref{15} and 
by point (i) of \lref{12}. Therefore $a^2 = 1$, and $a^x = [z,x\m,y]\m
[y,x\m,z] = a\m [z,x^2,y]\m][y,x^2,z] = a\m = a$. Similarly, $a^y = a$
and $a^z = a$. Finally, $[x,y,z][x,z,y]\m = [x,y,z][y,z,x]\m a = a
[x,y,z][x,y,z]^x = a [x^2,y,z] = as_x$.
\edem

We finish this section by an easy (but handy) observation:
\begin{lem}\label{17}
Let $Q$ be a Buchsteiner loop with elements $x$, $y$ and $z$. If
$[x,y,z]=1$, then $[y,z,x] = [z,x,y] = 1$ as well.
\end{lem}
\dem Use equalities $[x,y,z]^x = [y,z,x]\m$ and $[x,y,z]^{xy} =
[z,x,y]$.
\edem

\section{Central elements and an odd order proposition} \label{2}

As we have already hinted in the introduction, Buchsteiner loops are closely
connected to conjugacy closedness. By \cite{pac}, a loop $Q$ is conjugacy
closed if and only if $Q/N$ is an abelian group and all associators
are invariant under every permutation of their arguments. If we assume that
$Q/N$ is a boolean group, then the condition for Buchsteiner identity
is a natural weakening of the condition for the conjugacy closedness. Indeed,
if $Q/N$ is a boolean group, then $Q$ is a Buchsteiner loop if and only
if $[x,y,z] = [y,z,x]$ for all $x,y,z \in Q$, i.~e.~if the associators
are invariant under the cyclic shifts, by \cite{bcc} (see also \lref{46}).

By \cite{ccc}, if $Q$ is a Buchsteiner loop, then $Q/Z(Q)$ is a conjugacy
closed loop.

\begin{prop}\label{21}
Let $Q$ be a Buchsteiner loop with elements $x$, $y$ and $z$. Then
$[x,y,z][x,z,y]\m = [x,y,z][y,x,z]\m$, $[x,y,z][z,y,x]\m$, and
$[x^2,y,z] = [y,z,x^2]=[z,x^2,y]$ are central elements of
exponent $2$, and $[x,y,z]^{u^2} = [x,y,z]$ for each $u \in Q$.
\end{prop}
\dem Since the associators of a CC loop are invariant to permutations
of arguments,
there must be $[x,y,z]\equiv [x,z,y]$ and $[x,y,z]\equiv [z,y,x]
\bmod Z(Q)$, and so the initial claims follow from \pref{16}.
Since $Q/Z(Q)$ is a conjugacy closed Buchsteiner loop, each square
element belongs to the nucleus of $Q/Z(Q)$. Associators that involve
a nuclear element are trivial. Hence $[x^2,y,z] \equiv 1 \bmod Z(Q)$,
and so $[x^2,y,z] \in Z(Q)$. Furthermore, $1 = [x^4,y,z] = [x^2,y,z]^{x^2}
[x^2,y,z] = ([x^2,y,z])^2$, $[z,x^2,y]^z = [x^2,y,z]\m = [x^2,y,z]$ and
$[y,z,x^2]^y = [z,x^2,y]$. The last equality follows from expressing
$[xu^2,y,z]$ both as $[x,y,z]^{u^2}[u^2,y,z]$ and as
$[x,y,z][u^2,y,z]^x = [x,y,z][u^,y,z]$.
\edem

Note that \lref{11} is a special case of \pref{21}. Methods of \secref{1}
suffice to prove \pref{21} in Buchsteiner loops that are generated by
three elements, but it is an open question if these methods can be used
to prove \pref{21} in the full generality. To formalize this problem
consider a first order theory that involves a group $G \cong Q/A$, $A = A(Q)$,
that
acts upon a group $N$, and a ternary mapping $[-,-,-]: G^3 \to A \le Z(N)$.
In this theory we assume  that $[x,y,zu] =
[x,y,z]^u[x,y,u]$ and $[x,y,z]^x = [y,z,x]\m$
for all $x,y,z,u \in G$, and that $N/A$ can be identified
with a subgroup $H \le G$ in such
a way that $[-,-,-]$ depends only upon classes modulo $H$,
and $G/H$ is an abelian group of exponent four.

The associator calculus developed in \cite{buch} (which is an earlier
paper than \cite{ccc}) can be formulated within such a theory, and this
is also true for results of \secref{1}. The main results of this paper
are independent of \pref{21} since for them it suffices to know the
statement only for 3-generated groups.

For a commutative group  $G$ denote by $O(G)$ the subgroup consisting
of all elements of an odd order.
If a loop $Q$ contains a normal subloop $H$ which is a group, then
every characteristic subgroup of $H$ is clearly also a normal subloop of $Q$.
In particular, if $A(Q)$ is abelian, then $O(A(Q)) \unlhd Q$.

We have already mentioned that if a loop $Q$ is modulo the nucleus an
abelian group, then it is conjugacy closed if and only if each
$[x,y,z]$ does not depend on the order of the arguments. To verify
the latter property it suffices to show $[x,y,z] = [y,x,z]$
and $[x,y,z] = [x,z,y]$, for all $x,y,z \in Q$.

\begin{prop}\label{22}
Let $Q$ be a Buchsteiner loop that is not conjugacy closed. Then
neither $Q/O(A(Q))$ is conjugacy closed.
\end{prop}
\dem  We have $a=[x,y,z][y,x,z]\m =
[z,x,y][z,y,x]\m$, by \pref{16}, and so our assumption implies the
existence of $x,y,z \in Q$ such that $a \ne 1$. But then $a$ is an involution,
again by \pref{16}, and hence $a \notin O(A(Q))$.
\edem

\begin{lem}\label{23}
Let $Q$ be a Buchsteiner loop generated by a set $X$. If
$[x,y,z] = [x,z,y]$ for all $x,y,z \in X$, then $Q$ is a conjugacy closed
loop.
\end{lem}
\dem We need to prove $[t_2,t_1,t_3] = [t_1,t_2,t_3] = [t_1,t_3,t_2]$
for all $t_i \in Q$, $1\le i \le 3$. Since $Q$ is assumed to be 
a Buchsteiner loop, it suffices to prove only the latter identity.
Indeed, If $[t_3,t_1,t_2] = [t_3,t_2,t_1]$, then $[t_1,t_2,t_3] = 
([t_3,t_1,t_2]\m)^{t_3} = ([(t_3,t_2,t_1]\m)^{t_3} = [t_2,t_1,t_3]$,
The elements $t_i \in Q$ can be
regarded as terms in an abelian group of exponent $4$, for which $X$
is a set of generators. Each $t_i$ has thus a length $|t_i| \le 3|X|$,
and we can proceed by induction along $s=\sum |t_i|$. The case
$s=3$ is a consequence of our starting assumption. Let us have
$s\ge 4$. Then one of $t_i$, say $t_2$ is of the form $uv$. Using
the induction assumption we get $[t_1,uv,t_3] = [t_1,u,t_3]^v[t_1,v,t_3]
= [t_1,t_3,u]^v[t_1,t_3,v] = [t_1,t_3,uv]$.
\edem

\begin{cor}\label{24}
Let $Q$ be a Buchsteiner loop generated by $x$ and $y$. Suppose that
$Q$ is not conjugacy closed. Then $[x,x,y] \ne [x,y,x]$ or
$[y,y,x] \ne [y,x,y]$.
\end{cor}
\dem If $[x,x,y] = [x,y,x]$ and $[y,y,x] = [y,x,y]$, then 
$Q$ is conjugacy closed, by \lref{23} with $X = \{x,y\}$.
\edem

\begin{cor}\label{25}
Let $Q$ be a Buchsteiner loop generated by a set $X$. Let $Q_1$ be
the subloop generated by $X\setminus N$. If $Q_1$ is conjugacy
closed, then $Q$ is conjugacy closed as well.
\end{cor}
\dem This follows from \lref{23} too, since $[x,y,z] = 1 = [x,z,y]$
whenever $N \cap \{x,y,z\} \ne \emptyset$.
\edem 

\section{Loops that are not boolean modulo the nuclues}\label{3}

\begin{lem}\label{31}
Let $Q$ be a Buchsteiner loop with elements $x$, $y$ and $z$.
If $y^2 \in N(Q)$ and $z^2\in N(Q)$,
then $[x^2,y,z] = [y,z,x^2] = [z,x^2,y] =1$.
\end{lem}
\dem 
Both $[y^2,z,x]$ and $[z^2,y,x]$ are trivial, by our assumptions.
From \pref{13} we get $1 = [x^2,y,z][y^2,z,x][z^2,y,x] = [x^2,y,z]$.
The cyclic shifts of the latter associator are trivial by \lref{17}.
\edem

\begin{lem}\label{32}
Let $Q$ be a Buchsteiner loop such that $Q/N$ contains exactly one
nontrivial square element $x^2N$. Then there exists $y \in Q$ such
that $[x^2,x,y] \ne 1$, $y \notin xN$ and $y^2 \in N$.
\end{lem}
\dem 
If $[x^2,x,y] \ne 1$, then  $y \equiv x \bmod N$ since
$[x^2,x,x] = 1$, by Lemmas~\ref{11} and~\ref{14}, and
there must be $y^2 \in N$, by the assumption of the unique
square of $Q/N$. Therefore it suffices to find $y \in Q$
with $[x^2,x,y] \ne 1$.

The element $x^2$ does not belong to $N$, and hence $[x^2,y,z] \ne 1$
for some $y,z \in Q$. However, that means that at least one of $y^2$ and
$z^2$ does not belong to $N$, by \lref{31}. We can assume $z^2 \notin N$,
since $[x^2,y,z] = [x^2,z,y]$, by \pref{16}. If also $y^2 \notin N$, 
then $y \equiv xu\bmod N$ for some $u \in Q$ with $u^2 =1$. 
In such a case $[x^2,y,z] = [x^2,xu,z]
= [x^2,x,z]^u[x^2,u,z]$, and we are done if $[x^2,x,z] \ne 1$. Let us have
$[x^2,x,z] =1$. Then we are back to the case $[x^2,y,z] \ne 1$, but
now we can assume that
$y^2 \in N$. We know that $z \equiv xv$ for some $v \in Q$ with $v^2 \in N$
since $z^2 \notin N$, and there is a unique nontrivial square in 
$Q/N$. We have $[x^2,y,v]=1$, by \lref{31}, and so 
$1 \ne [x^2,y,z] = [x^2,y,x] = [x^2,x,y]$. The last equality
foolows from \lref{11}.
\edem

\begin{prop}\label{33}
Let $Q$ be a Buchsteiner loop generated by elements $x$ and $y$.
Suppose that $[x,x,y] \ne [x,y,x]$. Then $[x^2,x,y] \ne 1$,
$|Q/N| \ge 8$ and $|A(Q)| \ge 8$.
\end{prop}
\dem The equalities established in \lref{14} will be used
freely throughout the proof. Put $a =[x^2,x,y] = [x,x,y][x,y,x]\m$. 
We assume that $a\ne 1$,
and therefore $|Q/N|$ is not of exponent two. $|Q/N|$ cannot be cyclic,
since $[x^2,x,x] =1$, and thus $|Q/N| \ge 8$.
Now, $a^2 = 1$, by \lref{11}, and so $Q/O(A(Q))$ satisfies the hypothesis
(cf.~\pref{22}). We can hence assume that $A(Q)$
is a $2$-group. Set $u = [x,y,x]$. Then $u \ne a$ since $1\ne[x,x,y]\m
= u^x$. If $|u|$, the order of $u$, is greater than four, then
$|A(Q)| \ge 8$. Assume $|u| = 4$. If $a \ne u^2$, then $a$ and $u$
generate a subgroup of order $8$. Assume $a = u^2$. Then $u^x =
[x,x,y]\m = (ua)\m = u\m a = u a u^{-2} = u$, and so to prove $|A(Q)|
> 4$ it suffices to find an element $m \in A(Q)$ with $m^x \ne m$.
Set $m = [x,y,x][y,x,y]\m$. Then $m^x = [x,x,y]\m [y,x,y] = a m\m$,
and $m\m = m$ as $[x,y,x]^2 = [y,x,y]^2$.

It remains to consider the case when $u$ is an involution. In such a case
$u^y = u$. The element $a$ is central, by \lref{11}, and so to show
$|A(Q)| > 4$ it suffices to find $s \in A(Q)$ with $s^y \ne s$. Set
$s = [x,x,x]$. Then $s^y = s[x,x,y]\m[x,y,x]\m$, which equals $sa$ since
both $[x,x,y]$ and $[x,y,x]$ are assumed to be involutions.
\edem

The above proof can be seen as a starting point for constructing
Buchsteiner loops of order $64$ that are not boolean modulo the 
nucleus. As we prove below, $64$ is the least order for such 
a loop. An example was constructed in \cite{buch}, and one can
hope that all such loops of $64$ will be classified in future.
The overlap of \pref{33} and the ensuing \lref{34} should
be understood as justified by this intention.

\begin{lem}\label{34}
Let $x$, $y$ and $z$ be such elements of a Buchsteiner loop $Q$ that
satisfy $[x^2,y,z] \ne 1$. Then $|A(Q)| \ge 8$.
\end {lem}
\dem Set $s = [x^2,y,z]$. This is a central involution, and hence
$Q/O(A(Q))$ satisfies the hypothesis, and we can 
assume that $A(Q)$ is a $2$-group. Set $u = [x,y\m,z]$ and note
that $s = u^xu$, by \pref{13}. Therefore $u \ne 1$ and $u \ne s$.
If $|u| \ge 8$, then we are done. Let us have $s = u^2$. Then $u^x = u$,
and so $|A(Q)| \ge 8$ if we find $v \in A(Q)$ with
$v^x \ne v$.  Suppose that no such $v$ exists.  By considering
the formula of \lref{15} for $[y,z,y]^x$ and $[z,y,z]^x$,
we get $[y^2,z,x] = [z^2,x,y] = u^2$, by point (i) of \lref{12}
(and by \pref{16}). Of course, $[x^2,y,z] = [x,y,z]^2 = u^2$ as well.
From \pref{13} we now obtain $1 = (u^2)^3 = u^2$, a contradiction.

It remains to consider the situation when $u^2 = 1$
and $|A(Q)| = 4$. Then $u^x = us$,
and each of $u^y, u^z \in \{u,us\}$ since the element $s$
is central. We cannot have both $u^y$
and $u^z$ equal to $us$ because then $u^xu^yu^z = us \ne u$,
and that contradicts \pref{13}.

Assume $u^y = u$.
Then it suffices to find an element $v\in A(Q)$ with $v^y \ne v$.
By \lref{15}, $[x,z,x]^y = [x,z,x]s$ since $[y,z,x]^2 = 1$
by \lref{12}, and so
one can set $v= [x,z,x]$. The case $u^z = u$ is nearly the same.
\edem

\begin{cor}\label{35}
Let $Q$ be a Buchsteiner loop such that $Q/N(Q)$ is not a boolean
group. Then $|A(Q)| \ge 8$ and $|Q:N(Q)| \ge 8$. If $|Q| = 64$,
the $|A(Q)| = 8$ and $Q/N(Q) \cong C_4 \times C_2$.
\end{cor}
\dem Since $Q/N$ is not boolean, there must exist elements $x,y,z \in Q$
that satisfy the hypothesis of \lref{34}. Hence $|A(Q)| \ge 8$,
and so $Q/N$ has to be generated by at most two elements if $|Q| \le 64$.
In such a case we can assume that $Q$ is generated by two elements,
by \cref{25}, and we can also assume that $[x,x,y] \ne [x,y,x]$,
by \cref{24}. The inequality $|Q:N|$ now follows from \pref{33}
(or directly from \lref{32}).
\edem

For future references we also record this in a somewhat less explicit way:
\begin{cor}\label{36}
Let $Q$ be a Buchsteiner loop such that $Q/N$ is generated by less
than three elements. If $Q$ is not conjugacy closed, then $|Q/N| \ge 8$
and $|A(Q)| \ge 8$.
\end{cor}

\section{Commutator calculus and loops of small order}\label{4}

Let $x$ and $y$ be elements of a loop $Q$. The commutator $[x,y]$
is defined by $yx[x,y] = xy$. Assume that $N = N(Q) \unlhd Q$
and that $Q/N$ is an abelian group. Then $xy = yx[x,y]\m$ and so
$[y,x]=[x,y]\m$, as in groups. Furthermore, if $Q/N$ is an abelian
group, then one can connect associators and commutators by the formula
\begin{equation*}
 [xy,z] = [x,z]^y[y,z][x,z,y]\m [x,y,z][z,x,y].
\end{equation*}
The proof is not difficult, and can be found, e.~g., in \cite{buch}.

\begin{lem}\label{41}
Let $Q$ be a Buchsteiner loop with elements $x$, $y$ and $z$.
Set $m = [z,x,y][z,y,x]\m$. Then $m^2 =1$, $m \in Z(Q)$,
\begin{equation*}
[xy,z] = [x,z]^y[y,z][y,z,x]m
\text{ and }
[yx,z] = [y,z]^x[x,z][x,z,y]m.
\end{equation*}
\end{lem}
\dem By \pref{16} we can replace in the above formula the product
$[x,z,y]\m[x,y,z]$ with the product $[y,z,x][z,y,x]\m$.
That gives the required expression of $[xy,z]$,
and the expression of $[yx,z]$ uses the fact
that $m$ is a (central) element of exponent two, by \pref{21}. 
\edem

We shall apply \lref{41} to various situations, starting
with cases that naturally imply $[y,z,x] = [x,z,y]$.
The following observation be useful.

\begin{lem}\label{42}
Let $x$, $y$ and $z$ be elements of a Buchsteiner loop $Q$ such that
$[x,y,z]$ is centralized by each of $x$, $y$ and $z$. Then
$[x,y,z] = [y,z,x] = [z,x,y]$ is of exponent two, and
$[x^2,y,z]=[y^2,z,x] = [z^2,x,y] = 1$.
\end{lem}
\dem Use the notation of \pref{13}. We see that $u = [x,y\m,z]$
satisfies both $u^4 = 1$ and $u^6 = 1$ since we assume $u^x = u^y = u^z = u$,
and $u^2 = [x,y,z]^2 = s_y$ is of exponent two, by \lref{11}.
Hence $u^2=1$, and elements $s_x$, $s_y$ and $s_z$ are equal to $1$.
\edem

\begin{prop}\label{43}
Let $Q$ be a Buchsteiner loop with elements $x$, $y$ and $z$
such that all elements $[x,y]$, $[y,z]$ and $[x,z]$ are central.
Then $[x,y,z] = [y,x,z]$.
\end{prop}
\dem First note that $[xy,z] = [yx,z]$ since $xy = yx[x,y]$
and we assume $[x,y] \in Z(Q)$. The rest follows from \lref{41}.
\edem

\begin{cor}\label{44}
Let $Q$ be a Buchsteiner loop of nilpotency
class two. Then $Q$ has to be conjugacy closed.
\end{cor}
\dem The assumptions of both \lref{42} and \pref{43} are satisfied by all
$x,y,z \in Q$, and so we see that the value of an associator does not depend
upon the order of its arguments.
\edem

\begin{cor} \label{45}
Let $Q$ be a Buchsteiner loop such that $N(Q) \le Z(Q)$.
Then $Q$ has to be conjugacy closed.
\end{cor}
\dem Such a loop is necessarily of nilpotency class at most two.
\edem

\begin{lem} \label{46}
Let $Q$ be a Buchsteiner loop such that $Q/N$ is a boolean group.
Then
$$[x,y,z] = [y,z,x] = [z,x,y] = [x,y\m,z] \text{ for all }
x,y,z \in Q.$$
Furthermore, $[x,y,z]^x = [x,y,z]^y = [y,z,x]^z = [x,y,z]\m$.
\end{lem}
\dem
This follows directly from \lref{12} and \pref{13}.
\edem

\begin{prop}\label{47}
Let $Q$ be a Buchsteiner loop such that $|A(Q)| > 2$ and
$Q/N$ is a boolean group. If $Q$ is not conjugacy closed,
then $|Q|\ge 64$.
\end{prop}
\dem Throughout the proof we shall be assuming that 
$Q$ is not conjugacy closed. Thus
$|Q:N(Q)| \ge 8$, by \cref{25}, \cref{24} and \lref{14}. 
There cannot be $|N(Q)|=2$, since otherwise $N(Q)$
would be central, and \cref{45} would apply. We also
know that $|A(Q)$ is even, by \pref{22}. Choose 
$x,y,z \in Q$ such that $[x,y,z] \ne 1$ and $[x,z,y]
\ne [x,y,z]$, and denote by $Q_1$ the loop generated
by $x$, $y$ and $z$. If $A(Q_1)$ has only two elements,
then there must be $|Q:N(Q)| \ge 16$, and so $|Q| \ge 4\cdot 16
= 64$. We can hence assume $Q = Q_1$.

Our goal is to show that there must be $|N(Q)| \ge 8$.
Assume the contrary. If $|N(Q)| = 2$, then $N(Q) \le Z(Q)$,
and $Q$ is a CC loop, by \cref{45}. If $A(Q) = N(Q)$ is
of order $6$, then we obtain the same kind of contradiction,
by \pref{22}, and so $N(Q)=A(Q)$ has to consist of four elements,
and not all of them can be central.

Assume first that $A(Q)$ is a boolean group. To obtain
a contradiction, we shall show that $A(Q) \le Z(Q)$.
For that it suffices to prove $[u,v,w] \in Z(Q)$ for all possibilites
when $u,v,w \in \{x,y,z\}$, since the further cases follow from the associator
multiplicative formula. Now, if $\{u,v,w\} = \{x,y,z\}$, then
$[u,v,w] \in Z(Q)$ by \lref{46}. Furthermore, $[u,v,u]^w =
[u,v,u]$, by \lref{15}, and the rest follows from
\lref{14} in a clear way.

Let now $N(Q) =A(Q)$ be a cyclic group of order four. Denote by $b$
be the only nontrivial central element of $Q$. If $v\in Q$,
then $v^2 \in N(Q)$. If $v^2 \in Z(Q)$, then $[v,v,v]=[v^2,v]
= 1$. Consider an element $v \in Q$ with $[v,v,v] \ne 1$. Then
$u = v^2$ has to generate $N(Q)$, and $[v,v,v] = [v^2,v]$
is equal to $u\m u^v$. Thus $v$ has to induce the (only admissible)
nontrivial automorphism of $N(Q)$, and so $u^v = u\m$. That
means $[v,v,v] = b$, and so $[v,v,v] \in Z(Q) = \{1,b\}$ for
all $v \in Q$.

Consider elements $v,w \in Q$. We have $[v,v,v] = [v,v,v]^w$,
and the latter element is equal to $[v,v,v][v,w,v]^2$, by 
\lref{14}.  Thus $[v,w,v]^2 = 1$ for all
$v,w \in Q$, which in our situation means $[v,w,v] \in Z(Q)$.

To get a contradiction we shall prove now that $[u,v,w] \in Z(Q)$
for all $u,v,w \in Q$. This follows from \lref{15} since from
that lemma we see that $[u,v,w]^2 = 1$ for all $u,v,w \in Q$.
\edem

\begin{prop}\label{48}
Let $Q$ be a Buchsteiner loop of order less than $64$ that is not
conjugacy closed. Then $|Q| = 32$, $Q/N$ is elementary abelian
of order $8$, and $Z(Q) = A(Q)$ is of order $2$. The group $Q/Z(Q)$
is a nonabelian group of order $16$.
\end{prop}
\dem
From \cref{36} we know that $Q/N$ has to be of order at least $8$,
and from \cref{35} we know that it has to be elementary abelian.
Furthermore $|A(Q)| =2$, by \pref{47}. Thus $A(Q) \le Z(Q)$,
and there cannot be $A(Q) = N(Q)$, by \cref{45}. This means
that $|Q:N| = 8$ and $|N| = 4$. From \cref{45} we also see that
$Z(Q)$ has to coincide with $N(Q)$. Finally, $Q/Z(Q)$ cannot be
abelian, by \cref{44}.
\edem

\section{The doubling construction}\label{5}

The purpose of this section is to describe a construction based upon 
a Buchsteiner loop $Q$ that produces a Buchsteiner loop $P$ that
contains the loop $Q$ as a subloop of index two. It may happen
that $Q$ is a CC loop, while $P$ is not, and in the next section
we shall see that all proper Buchsteiner loops of order 32 can be
obtained in this way.

\begin{prop}\label{51}
Let $P$ be a Buchsteiner loop a with a normal subloop $Q$, where $|P:Q| = 2$.
Let $z\in P \setminus Q$ be an element such that $d = z^2 \in N(P)$
and such that $q(u)=[z,u]$  belongs to $Z(Q)$ and is of exponent two 
for all $u \in Q$.
Then:
\begin{enumerate}
\item[(i)] $[u,v,z] = [v,z,u] = [z,u,v] = q(vu)q(u)q(v)\in Z(P)$ for all $u,v \in Q$;
\item[(ii)] $[z,z,z] = q(d) \in Z(P)$;
\item[(iii)] $[u,z,z] = [z,u,z]=[z,z,u] =[d,u]=[u,d] \in Z(P)$.
\end{enumerate}
\end{prop}
\dem We assume $[z,u]^2 = 1$, and hence $q(u) = [u,z] = [z,u]$ for all $u\in Q$.
Consider $u,v \in Q$. Then $[uv,z] = [u,z]^v[v,z][u,v,z]m$, 
$m=[v,z,u][v,u,z]\m \in Z(P)$ and $m^2 = 1$, by \lref{41}.
The element $[u,v,z]$ can be thus expressed as a product of elements from
$Z(Q)$ that are of exponent two. Therefore $[u,v,z] = [u,v,z]^u = [v,z,u] = [v,z,u]^v
= [z,u,v]$, and we get $[uv,z] = [u,z][v,z][v,u,z]$. We also have $[u,v,z]^z = [z,u,v]\m
= [u,v,z]\m = [u,v,z]$, and thus $[u,v,z] \in Z(P)$.

For point (ii) it suffices to note that $[z,z,z]=[z^2,z] = [d,z] = 
q(d)\in Z(Q)$, and that $[z,z,z]^z = [z,z,z]\m=[z,z,z]$.

For each $u \in Q$, $[z,u,z] = [z,z,u]$, by \lref{14}, since $d= z^2\in N(P)$.
If $v \in Q$, then $[z,u,z]^v = [z,u,z][z^2,u,v][z,u,v]^{-2} = [z,u,z]$,
by \lref{15} and  by point (i) of this proof. This also gives
$[z,u,z]^2 = 1$ since $[z,u,z]^u = [z,u,z]\m$, by \lref{14}.
Furthermore, $[z,uz,z] = [z,z,z]^u[z,u,z]=[z,z,z][z,u,z]$,
and $[z,uz,z] = [z,z,z][z,u,z]^z$. Therefore $[z,u,z]^z = [z,u,z]$,
and so $[z,u,z] \in Z(P)$.
Finally,  \lref{41} yields $[d,u] = [z^2,u] = [z,u]^2[z,u,z] = [z,u,z]$ .
\edem

For the next few statements we shall assume that $P$ is as in \pref{51}.
The associator multiplicative formulas immediately imply:

\begin{cor}\label{52}
Let $A$ be a subloop of $P$ generated by all associators $[\alpha,\beta,\gamma]$
such that $z \in \{\alpha,\beta,\gamma\}$. Then $A$ is a boolean group that
is contained in $Z(P)$.
\end{cor}
\begin{cor}\label{53}
Assume $u_i \in Q$ and $\eps_i \in \{0,1\}$, $1 \le i \le 3$. Then
\begin{multline*}
[u_1z^{\eps_1},u_2z^{\eps_2},u_3z^{\eps_3}] =
[u_1,u_2,u_3][z,u_2,u_3]^{\eps_1}[u_1,z,u_3]^{\eps_2}
[u_1,u_2,z]^{\eps_3} \\
[u_1,z_2,z_3]^{\eps_2\eps_3} [z_1,u_2,z_3]^{\eps_1\eps_3}
[z_1,z_2,u_3]^{\eps_1\eps_2}[z_1,z_2,z_3]^{\eps_1\eps_2\eps_3}.
\end{multline*}
Furthermore, $[u_1,u_2,u_3]^z = [u_1,u_2,u_3]$.
\end{cor}
\dem Only the last equality requires a proof. We have
$[u_1,u_2,u_3]^z[z,u_2,u_3] = [u_1z,u_2,u_3] = [u_1,u_2,u_3][z,u_2,u_3]^{u_1}
= [u_1,u_2,u_3][z,u_2,u_3]$, by \cref{52}.
\edem

\begin{lem}\label{54}
The loop $Q$ contains normal subloops $A\le S$ such that $A\le Z(Q)$,
$N(P)\cap Q \le S$, both $A$ and $Q/S$ are boolean groups, and there exist
mappings $q:Q\to A$ and $\vhi: Q \to A$ such that:
\begin{enumerate}
\item[(i)] $q(a) = 1$ for all $a \in A(Q)$;
\item[(ii)] $q(u) = q(u') = [z,u]$ whenever $u \equiv u' \bmod A$,
for all $u,u'\in Q$;
\item[(iii)] $q(du) = q(d)q(u)$ for all $u \in Q$;
\item[(iv)] the mapping $g(u,v) = q(vu)q(u)q(v)$ induces, for all $u,v \in Q$,
a group homomorphism $Q/S \to A$ whenever one of the coordinates is fixed; and
\item[(v)] the mapping $\vhi(u) = [u,d]$ induces a group
homomorphism $Q/S \to A$.
\end{enumerate}
\end{lem}
\dem Let $A$ be defined as in \cref{52}. We have $A \le Z(P)\cap Q \le Z(Q)$,
and $A$ is a boolean group. If $a \in A(Q)$, then $q(a) = [z,a] = [a,z] =
a\m a^z$ since $a\in A(P)\le N(P)$, and $a^z = a$, by \cref{53}. This
proves point (i).

Point (ii) is clear since $[z,u] = [z,ua]$ for all $u\in Q$ and  $a \in A$ 
as $A \le Z(P)$, by \cref{52}.

We assume $z^2 \in N(P)$, and so $[du,z] = [d,z]^u[u,z]=[d,z][u,z]$, 
by \lref{41}. That gives (iii).

Now, $g(u,v) = [u,v,z]$ for all $u,v \in Q$, by point (i) of \pref{51}.
The values of $g(u,v) = [u,v,z]\in Z(P)$ depend only upon classes of $u$ and $v$
modulo $N(P)$, and thus $g(u,v) = g(u',v')$ when $u\equiv u'$
and $v\equiv v' \bmod N(Q)$. The 
multiplicative associator formula immediately implies 
$g(u_1u_2,v) = g(u_1,v)g(u_2,v)$ and $g(u,v_1v_2) = g(u,v_1)g(u,v_2)$,
for all $u,v,u_1,u_2,v_1,v_2 \in Q$. Similarly, $\vhi(uv)
= \vhi(u)\vhi(v)$ for all $u,v \in Q$ since $\vhi(-) = [z,-,z]$,
by point (iii) of \pref{51}.

We thus have homomorphisms $[z,-,z]$, $[-,u,z]$ and
$[u,-,z]$ that map $Q$ into $A$, where $u$ runs through $Q$. 
Since $A$ is a boolean group, 
$Q$ is a boolean group modulo the kernel of each homomorphism.
Therefore it is boolean modulo the intersection of all kernels,
and this intersection gives a subloop $S\ge N(P)\cap Q \ge N(Q)$ that is
required by our statement. Points (iv) and (v) thus follow from
the earlier computations and from the latter fact.
\edem

\begin{lem}\label{55}
The following equalities hold for all $u,v \in Q$:
\begin{gather*}
zu\cdot v = z\cdot uvg(u,v),\  u\cdot zv = z\cdot uvq(u)g(u,v)g(v,u), \\
\text{and \,}zu\cdot zv = duv\vhi(v)q(u)g(v,u).
\end{gather*}
\end{lem}
\dem Indeed, one can compute easily that 
\begin{align*}
zu \cdot v &= z\cdot uv[z,u,v] = z\cdot uv g(u,v), \\
u \cdot zv &= uz\cdot v [u,z,v] = zu[u,z]\cdot v [u,z,v] = z\cdot uv [u,z]
[z,u,v][u,z,v] \\
&= z\cdot uv q(u)g(u,v)g(v,u), \text{ and} \\
zu \cdot zv &= [z,u]uz\cdot zv = u(z \cdot zv)[z,u][u,z,zv] = u(z^2v)[z,z,v][z,u][u,z,zv] \\
&= [u,d]duv[d,v] q(u)[u,z,v][u,d] = duv\vhi(v)q(u)g(v,u).
\end{align*}
\edem

Our aim now is to show that the properties of the above loop $P$ can
be used for a construction based on $Q$, $d$ and $q$.

Suppose that $Q$ is a loop with normal subloops $A$ and $S$.
Suppose that $d$ is an element of $N(Q)$, and that $q: Q \to A$ a mapping.
Put $\vhi(u) = [d,u]$ for all $u \in Q$, and
$g(u,v) = q(vu)q(u)q(v)$ for all $u,v \in Q$.
Assume that
\begin{enumerate}
\item[(1)] both $A$ and $Q/S$ are boolean groups, and
$A \le S\cap Z(Q)$;
\item[(2)] $q(ua) = q(u)$ for all $u \in Q$ and $a \in A$;
\item[(3)] $g(u,vw) = g(u,v)g(u,w)$ and $g(vw,u) = g(v,u)g(w,u)$
for all $u,v,w \in Q$;
\item[(4)] $g(us,vt) = g(u,v)$ for all $u,v \in Q$ and $s,t \in S$;
\item[(5)] $\vhi(uv) = \vhi(u)\vhi(v)$ for all $u,v \in Q$;
\item[(6)] $\vhi(us) = \vhi(u)$ for all $u \in Q$ and $s \in S$; and
\item[(7)] $q(du) = q(d)q(u)$ for all $u \in Q$.
\end{enumerate}
Define a loop $P(*) = Q[d,q,z]$ on $Q\cup zQ$, $z \notin Q$, by
\begin{alignat*}{2}
u * v &= uv, & u*zv &= z\cdot uvq(u)g(u,v)g(v,u), \\
zu * v &= zuvg(u,v), &\quad \text{ and\quad } zu * zv &= duv\vhi(v)q(u)g(v,u),
\end{alignat*}
for all $u,v \in Q$.

The notation $Q[d,q,z]$ does not
carry an identification of subloops $A$ and $S$. This is not needed, indeed,
since $A$ can be replaced by the (central boolean) subgroup generated by
all $q(u)$, and $S$ can be replaced by the set of all $x \in Q$ such
that $g(u,x) = g(x,u) = 1$ for all $x \in Q$.

The following three lemmas are stated under the assumption that $P$ and 
$Q$ are as in the above construction.

\begin{lem}\label{56}
Each element of $A$ belongs to $Z(P)$.
\end{lem}
\dem Consider $u,v \in Q$ and $a \in A$. Then $(u*v)*a
= uva = u*(v*a)$, $zu * a = zua$, $zu * (v*a) =
z\cdot uvag(u,va) = z \cdot uvg(u,v)a =(zu*v)*a$,
$u * (zv *a) = z \cdot uvaq(u)g(u,va)g(va,u) =
z\cdot uvq(u)g(u,v)g(v,u)a = (u*zv)*a$ and
$zu*(zv *a) = duva\vhi(va)q(u)g(va,u) = duv\vhi(v)
q(u)g(v,u)a = (zu * zv) *a$, which means that $a$
belongs to the right  nucleus. Clearly, $zu * a =
zua = a * zu$, and so it remains to show that $a$
also belongs to the left nucleus. This follows
from $a*(zu * v) = z \cdot uvg(u,v)a = z \cdot uavg(ua,v)
= zua * v = (a*zu)*v$, $a * (u * zv) = (u*zv)*a = ua *zv
= (a*u)*zv$ and $a * (zu * zv) = (zu * zv) * a = zua * zv
= (a * zu) * zv$.
\edem

\begin{lem}\label{57}
Let $u$, $v$ and $w$ be elements of $Q$. Then
\begin{enumerate}
\item[(i)] $[zu,v,w] = [u,v,w]g(v,w)$,
\item[(ii)] $[u,zv,w] = [v,w,u]g(w,u)$,
\item[(iii)] $[u,v,zw] = [u,v,w]g(u,v)$,
\item[(iv)] $[zu,zv,w] = [u,v,w]\vhi(w)g(w,u)g(v,w)$;
\item[(v)] $[zu,v,zw] = [u,v,w]\vhi(v)g(v,w)g(u,v)$;
\item[(vi)] $[u,zv,zw] = [u,v,w]\vhi(u)g(w,u)g(u,v)$; and
\item[(vii)] $[zu,zv,zw] = [u,v,w]\vhi(d)\vhi(uvw)g(u,v)g(v,w)g(w,u)$.
\end{enumerate}
\end{lem}
\dem Our goal is to compute $[\alpha,\beta,\gamma]$, where
$\alpha \in \{u,zu\}$,
$\beta \in \{v,zv\}$ and $\gamma \in \{w,zw\}$. By using the definition of $*$
we shall in every case first express $(\alpha * \beta) * \gamma$ as
$z^\eps d^\eta(uv\cdot w) a(\alpha,\beta,\gamma)$, where $\eps,\eta \in\{0,1\}$
and $a = a(\alpha,\beta,\gamma) \in A$. Then we express
$\alpha * (\beta * \gamma)$
as $z^\eps d^\eta(u\cdot vw) b(\alpha,\beta,\gamma)$, where
$b = b(\alpha,\beta,\gamma)
\in A$ again. Now, $[\alpha,\beta,\gamma]$ should be equal to
$[u,v,w]c(\alpha,\beta,
\gamma)$, with $c = c(\alpha,\beta,\gamma) \in A$. To prove
$(\alpha *(\beta*\gamma))[u,v,w]c = (\alpha * \beta) *\gamma$ we hence need to
show that $(u\cdot vw)b[u,v,w]c = (uv\cdot w)a$, which amounts to $bc = a$,
which is the same as $abc = 1$.

In case (i) we get $(zu * v) * w = z(uv \cdot w)g(u,v)g(uv,w)$ (since
$g(uvg(u,v),w) = g(uv,w)$) and $zu * (v * w) = z(u\cdot vw)g(u,vw)$.
We have to verify that the product $g(u,v)g(uv,w)g(u,vw)g(v,w)$ vanishes,
and that clearly
follows from the equality
$g(uv,w)g(u,vw) = g(u,w)g(v,w)g(u,v)g(u,w) = g(v,w)g(u,v)$.

To get (ii) compute $(u * zv) * w = z(uv \cdot w)q(u)g(u,v)g(v,u)g(uv,w)$,
$u*(zv * w) = z \cdot(u\cdot vw)g(v,w)q(u)g(u,vw)g(vw,u)$ and
\begin{multline*}
q(u)g(u,v)g(v,u)g(uv,w)g(v,w)q(u)g(u,vw)g(vw,u)g(w,u) =\\
g(u,v)g(v,u)g(u,w)g(u,v)g(u,w)g(v,u) = 1.
\end{multline*}

For (iii) we get $(u*v)*zw = z\cdot (uv\cdot w)q(uv)g(uv,w)g(w,vu)$,
$u*(v*zw) = z \cdot (u\cdot vw) q(v)g(v,w)g(w,v) q(u)g(u,vw)g(vw,u)$
and
\begin{multline*}
q(uv)g(uv,w)g(w,vu)q(v)g(v,w)g(w,v) q(u)g(u,vw)g(vw,u)g(u,v) = \\
g(v,u)g(u,w)g(w,u)g(u,w)g(v,u)g(w,u) = 1.
\end{multline*}

To verify (iv) observe that $(zu * zv) * w = d(uv\cdot w)\vhi(v)
q(u)g(v,u)$, $zu * (zv *w) = d(u\cdot vw)\vhi(vw)q(u)g(vw,u)g(v,w)$
and
\begin{multline*}
\vhi(v)q(u)g(v,u)\vhi(vw)q(u)g(vw,u)g(v,w)\vhi(w)g(w,u)g(v,w) = \\
g(v,u)g(v,u)g(w,u)g(v,w)g(w,u)g(v,w) = 1.
\end{multline*}

Point (v) follows from $(zu * v)*zw = d(uv\cdot w)\vhi(w)q(uv)g(uv,w)
g(u,v)$, $zu * (v * zw) = d(u\cdot vw)\vhi(vw)q(u)g(vw,u)q(v)g(v,w)
g(w,v)$ and
\begin{multline*}
\vhi(w)q(uv)g(uv,w)g(u,v)\vhi(vw)q(u)g(vw,u)q(v)g(v,w)g(w,v)
\vhi(v)g(v,w)g(u,v) \\ =
g(v,u)g(u,w)g(v,w)g(u,v)g(v,u)g(w,u)g(w,v)g(u,v) = 1.
\end{multline*}

To get (vi) note that $(u * zv) * zw = d(uv\cdot w)
\vhi(w)q(uv)g(w,uv)q(u)g(u,v)g(v,u)$, $u *(zv * zw) =
d(u\cdot vw)\vhi(u)\vhi(w)q(v)g(w,v)$ and
\begin{multline*}
\vhi(w)q(uv)g(w,uv)q(u)g(u,v)g(v,u)\vhi(u)\vhi(w)q(v)g(w,v)
\vhi(u)g(w,u)g(u,v) \\ =
g(v,u)g(w,u)g(w,v)g(u,v)g(v,u)g(w,v)g(w,u)g(u,v) = 1.
\end{multline*}

Finally,
$(zu * zv) *zw = zd(uv\cdot w)
q(uv)g(uv,w)g(w,uv)\vhi(v)q(u)g(v,u)$, $zu * (zv * zw)=
zd (u \cdot vw)\vhi(u)\vhi(w)q(v)g(w,v)g(u,vw)$ and
\begin{multline*}
q(uv)g(uv,w)g(w,uv)\vhi(v)q(u)g(v,u)\vhi(u)\vhi(w)q(v)g(w,v)g(u,vw)
\\ \vhi(uvw)g(u,v)g(v,w)g(w,u)  =
g(v,u)g(u,w)g(v,w)g(w,u)g(w,v)
\\ g(v,u)g(w,v)g(u,v)g(u,w)
g(u,v)g(v,w)g(w,u) = 1.
\end{multline*}
\edem

\begin{lem}\label{58}
Suppose that $N(Q) \unlhd Q$ and that $Q/N(Q)$ an abelian
group. Then $A(P) \le N(P) \unlhd P$, with
$P/N(P)$ an abelian group. Furthermore, $N(Q) \cap S \le N(P)$.
\end{lem}
\dem If $x \in N(Q) \cap S$, then $[x,\alpha,\beta] =
[\alpha,x,\beta] = [\alpha,\beta,x] = 1$, for all $\alpha,\beta \in P$.
This follows directly from \lref{57}, by inspecting all possible
situations that are described by points (i)-(vi). Hence $N(Q) \cap S
\le N(P(*))$. From \lref{57} we also see that $A(P(*)) \le N(Q) \cap S$.
For the rest of the proof it suffices to find a commutative group $G(*)$
and a homomorphism $f:P(*) \to G(*)$ such that $N(Q) \cap S$ is
equal to the kernel of $f$.

Put $\bar Q = Q/(S\cap N(Q))$. Then $\bar Q$ is a commutative group,
as both $Q/S$ and $Q/N(Q)$ are assumed to be commutative groups.
Define now a loop $G(*)$ on $G = \bar Q \cup z\bar Q$ by
$z\bar u * \bar v = z\cdot \bar u \bar v$, $\bar u * z \bar v =
z\cdot \bar u \bar v$ and $z \bar u * z \bar v = \bar d \bar u \bar v$.
The operation $*$ is clearly commutative. To see that it is
associative one can use \lref{56}, with $\bar q$ and $\bar \vhi$
trivial, where $\bar A = 1$ and $\bar S = S/(S\cap N(Q))$.
The mapping $f$ is now defined by $f(u) = \bar u$ and $f(zu) = z\bar u$.
It is clear that this is a homomorphism $P(*) \to G(*)$ and
that $N(Q) \cap S$ is its kernel. \edem

\begin{lem}\label{59}
Suppose that $Q$ is a Buchsteiner loop such that $q([u,v,w])=1$
for all $u,v,w \in Q$.
Then $P$ is a Buchsteiner loop as well.
\end{lem}
\dem
The conditions of \lref{58} are satisfied and hence we know
that $A(P(*)) \le N(P(*)) \unlhd P(*)$. Therefore we only need to
prove that $[\alpha,\beta,\gamma]^\alpha =
[\beta,\gamma,\alpha]\m$, for all $\alpha,\beta,\gamma\in P$.
If $x \in Q$, then $z\bs (xz) = xq(x)$, by the definition of $P$. 
Thus $[u,v,w]^z = [u,v,w]$ for all $u,v,w \in Q$, by assumptions
of the lemma. The right hand sides
of all equalities in \lref{57} are hence invariant under the
action of $z$. This means that we need to verify $[\alpha,\beta,\gamma]^u
= [\beta,\gamma,\alpha]\m$ for all cases when $\alpha \in \{u,zu\}$,
$\beta \in \{v,zv\}$ and $\gamma \in \{w,zw\}$. Now, $[\alpha,\beta,
\gamma]^u = [u,v,w]^uc(\alpha,\beta,\gamma)$ for some
$c(\alpha,\beta,\gamma) \in A$, and $[\beta,\gamma,\alpha]\m
= [v,w,u]\m c(\beta,\gamma,\alpha)$. Since we assume $[u,v,w]^u
= [v,w,u]\m$, we have to show that $c(\beta,\gamma,\alpha)
= c(\alpha,\beta,\gamma)$, for all cases (i)-(vii) of \lref{56}.
Now, indeed $c(v,w,zu)=g(v,w)$, $c(zv,w,u) = g(w,u)$,
$c(v,zw,u) = g(u,v)$, $c(zv,w,zu) = \vhi(w)g(w,u)g(v,w)$,
$c(v,zw,zu) = \vhi(v)g(u,v)g(v,w)$, $c(zv,zw,u) = \vhi(u)
g(u,v)g(w,u)$, and the last case is clear since it is cyclically
invariant.
\edem

We are now ready for the final statements of this section.
\begin{prop} \label{510}
Let $Q$ be a Buchsteiner loop with normal subloops $A$ and $S$,
and with an element $d \in N(Q)$. Furthermore, let $q: Q \to 
A$ be a mapping such that $q(a) = 1$ for all $a \in A(Q)$, and let
$z$ be an element outside $Q$. If $d$ and $q$ satisfy conditions
(1)--(7), then $P = Q[d,q,z]$ is a Buchsteiner loop with $A\le Z(P)$,
$A(P) \le N(Q) \cap S \le N(P)$ and $N(Q) \cap S \le N(P)$, where
$P/N(Q)\cap S$ is an abelian group.

If $Q$ is a conjugacy closed loop, then $[\alpha,\beta,\gamma]
= [\beta,\gamma,\alpha]$ for all $\alpha,\beta,\gamma \in P$.
In such case $P$ is conjugacy closed if and only if $g(u,v)
= g(v,u)$ for all $u,v \in Q$.
\end{prop}
\dem Use Lemmas \ref{56}, \ref{58} and \ref{59}.

\begin{prop}\label{511}
Let $P$ be a Buchsteiner loop that contains a normal subloop $Q$,
$|P:Q|=2$, and an element $z \in P\setminus Q$ such that 
$d = z^2 \in N(P)$, and $[z,u]$ is a central element of $Q$,
$[z,u]^2 =1$, for all $u \in Q$. Set $q(u) = [z,u]$ for
all $u \in Q$. Then $P = Q[d,q,z]$.
\end{prop}
\dem This is just another expression of \lref{55}.
\edem

\section{Proper Buchsteiner loops of order 32}\label{6}

We shall first show that such loops really exist, by applying
the doubling construction of \secref{5} to the group $Q = G\times A$,
where $G$ is a group of quaternions and $A$ is a two-element group.
The (only) natural choice for $S$ is the subgroup $G' \times A$. 
The mapping $q:Q\to A$ has to depend only upon the elements
of $G$ (by condition (2)), and so we shall be looking for a mapping
$q: G\to \{0,1\}$ such that $g(u,v) = q(vu)+q(u)+q(v)$ yields a bilinear
mapping $G/G' \to \{0,1\}$. If $q$ is such a mapping, then we can always
set $d = 1$, and that gives a
a Buchsteiner loop $P$, by \pref{510}. However, the loop $P$ might be 
conjugacy closed. To avoid this case we need to make sure that
$g(-,-)$ is not symmetric (see \pref{510} again).

\begin{lem}\label{61}
Let $G$ be a group of quaternions generated by elements $x$, $y$ and $z$
such that $xy=z$, $yz = x$ and $zx = y$. Let $s = x^2 = y^2 = z^2$
be the only nontrivial square of $G$. Define $q:G \to \{0,1\}$ in
such a way that $q(u) = 1$ if and only if $u \in\{s,x,y,z\}$. Then
$G/G'$ is a vector space over $\{0,1\}$, and the mapping 
$g:G\times G\to \{0,1\}$, $(u,v) \to q(vu) + q(u) + q(v)$, induces
 a non-symmetric bilinear form on $G/G'$.
 \end {lem}
\dem We see that $q(us) = q(s)+ q(u)$ for all $u \in G$. The element
$s$ is central and so $g(u,v)$ clearly does not change if $u$
is replaced by $us$ or $v$ by $vs$. If $u$ and $v$ generate $G$,
then $s = [u,v]$, and $g(u,v) = g(v,u) + 1$.  For the
proof it therefore suffices to show that $g(u,vw) = g(u,v)+g(u,w)$,
where $u,v,w \in \{x,y,z\}$ and $vw \in \{s,x,y,z\}$. 
The case $v = w$ is clear, and so $v \ne w$
can be assumed. We can also assume $u = x$ because $\aut(G)$ acts
transitively upon $\{x,y,z\}$. Now $g(x,xy) = g(x,z) = q(zx) = 1 =
g(x,x) = g(x,x) + g(x,y)$, $g(x,zx) = 0 = g(x,x) + g(x,z)$, and
$g(x,yz) = g(x,x) = 1 = g(x,z) = g(x,y)+g(x,z)$.
\edem

\begin{cor}\label{62} There exists a proper Buchsteiner loop of 
order $32$.
\end{cor}
\dem Indeed, set $P = Q[q,1,z]$, with $Q = G\times A$, $G \cong Q_8$,
$A \cong C_2$, $z \notin Q$, and $q(ua) = 1$, where $u\in G$
and $a \in A$, if and only if $u\in \{1,x\m,y\m,yx\}$, for some
generators $x$ and $y$ of $G$.
\edem

\begin{thm}\label{63}
Let $P$ be a proper Buchsteiner loop of order $32$. Then
$1 < A(P) < N(P) < Z_2(P)$, $|Z_2(P)| = 8$, $A(P) = Z(P)$
and there
exists a unique power associative conjugacy closed subloop $Q$
of index two such that $QZ_2(P) = P$ and $Q\cap Z_2(P) = N(P)
= Z(Q)$. The group $Q/A(P)$ is noncommutative.
\end{thm}
\dem Set $A = A(P)$ and $N = N(P)$. We have $A=P$, $|A|=2$ and $|N|=4$,
by \pref{48}. Set also $C= Z_2(P)$. The subloop $C$ consists of all 
elements that are central modulo $A$. Group $H=P/A$ is nonabelian,
and group $P/N$ is elementary abelian of order $8$, again by \pref{48}.
The group $H$ thus contains a two-element subgroup modulo which
it is a vector space of dimension $3$, and the square mapping induces
a quadratic form of the vector space into this subgroup. The radical
of this quadratic form corresponds to $Z(H)$, and so $|Z(H)| = 4$.
The preimage of $Z(H)$ modulo $A$ is equal to $C$, the second centre
of $P$. If $x\in Q$ and $c\in C$, then $[x,c] \in Z=A$, and
$[xy,c] = [x,c][y,c][y,x,c]$, by \lref{41}. Furthermore,
$[c,u] = [u,c]$ and $[u,c]^2 = 1$, since $A$ has only two elements.
Clearly, $c^2 \in N$.

Consider the action of $H = P/A$ upon $N$, $n\mapsto x\bs (nx)$. 
This action has to be nontrivial since $A = Z$, and $|N:A|=2$. However,
each element of $H$ acts trivially upon $A$, and so the image
of the action contains exactly two permutations (the identity and 
the transposition of elements of $N\setminus A$). The kernel of this
action is hence a subgroup of $H$ that is of index two. The preimage
of the kernel in $H$ is a subloop $Q$, and this subloop satisfies
$Z(Q) \ge N$. Note that $P$ contains exactly one such subloop of index
two since each element of $Q$ acts trivially upon $N$. 

We shall be now establishing the properties of $Q$. It is clear that
$Q$ is conjugacy closed, by \pref{48}. Choose $x,y,z \in P$ so that
they form a basis modulo $N$, and $z \in C$. From \cref{25} we
see that these elements generate $Q$. The associator $[x,y,z]$ is central,
and hence invariant under cyclic shifts, by \lref{42}. Therefore
$[x,y,z] \ne [y,x,z]$, by \lref{23}, and hence $[xy,z]\ne [yx,z]$,
by the formula $[xy,c] = [x,c][y,c][y,x,c]$. Now, the same formula
gives $[xy,z] = [yx[x,y],z] = [yx,z][[x,y],z]$, and so we see that
$z$ acts nontrivially upon $[x,y] \in N$. That means that $z$ cannot
belong to $Q$. Thus $N = Q\cap C$. In fact, we have shown even more,
since for each $x \in Q\setminus N$ we can find $y \in Q$ such that
$x,y,z$ is a basis modulo $N$, and so for each $x \in Q\setminus N$
there exists $y \in Q$ with $[x,y] \ne 1$. Hence $Z(Q) = N$. The loop
$Q$ is power associative since $x^2 \in Z(Q)$ for all $x \in Q$.
\edem

\begin{cor}\label{64}
Each proper Buchsteiner loop of
order 32 can be obtained by the doubling construction.
\end{cor}

\section{Abelian inner mappings groups}\label{7}

We start by applying well known facts about inner mappings to
Buchsteiner loops.
\begin{lem}\label{71}
Let $Q$ be a Buchsteiner loop such that $A(Q) \le N(Q)$. Then
$$L(x,y)(z) = z [x,y,z]\m, \quad R(x,y)(z) = z[y,x,z], \text{\, and \ }
T_x\m(z) = z[z,x],$$
for all $x,y,z \in Q$.
\end{lem}
\dem Recall that $(x\cdot yz)[x,y,z] = xy \cdot z$. This means
$(x\cdot yz)= xy \cdot (z[x,y,z]\m)$, since $[x,y,z] \in N(Q)$,
and so $L(x,y)(z) = z[x,y,z]\m$. Now $(z\cdot yx)[z,y,x] =
zy \cdot x = (((z\cdot yx)[z,y,x])/(z\cdot yx))(z\cdot yx) =
[z,y,x]^{(zyx)\m}(z \cdot yx) = [z,y,x]^{x\m y\m z\m}(z\cdot yx) =
[y,x,z]^{z\m}(z \cdot yx) = (((z[y,x,z])/z)z)\cdot yx = z[y,x,z]
\cdot yx$. Hence $R(x,y)(z) = z[y,x,z]$. To prove $x\bs (zx) = z((xz)\bs (zx))$
it suffices to multiply the equality by $x$ on the left, and to use
the fact that $[z,x] = (xz)\bs(zx)$ belongs to the nucleus.
\edem

\begin{lem}\label{72}
Let $Q$ be a Buchsteiner loop. Then the set of all $L(x,y)$ and $R(x,y)$
generates an abelian group, and this group belongs to the center of
$\inn Q$ if and only if $A(Q) \le Z(Q)$.
\end{lem}
\dem Clearly, $R(x,y)R(u,v)(z) = R(x,y)(z[v,u,z])= z[v,u,z]
[y,x,z[v,u,z]] = z[v,u,z][y,x,z]$, and the other cases are similar
(in fact, their inspections is not needed when one takes in account
that $\ml_1 = \mr_1$, in every Buchsteiner loop $Q$). Now, $R(x,y)
T_u\m(z) = R(x,y)(z[z,u])= z[z,u][y,x,z[z,u]]=z[z,u][y,x,z]$, and
$T_u\m R(x,y)(z) = T_u\m (z[y,x,z]) = z[y,x,z][z[y,x,z],u]$.
Set $a=[y,x,z]$ and note that$a \in Z(N)$, and that $[za,u] =
[z,u][a,u]$, by \lref{41}. Hence $R(x,y)$ and $T_u$ commute for all $x,y,u
\in Q$ if and only if $[a,u] = 1$ for all $a \in A(Q)$. This is
the same as to say that $A(Q) \le Z(Q)$.
\edem

\begin{prop}\label{73}
Let $Q$ be a Buchsteiner loop with $A(Q) \le Z(Q)$. Then both $Q/N$
and $A(Q)$ are boolean groups and $[x,y,z] = [y,z,x]$ for all $x,y,z \in Q$.
If $A(Q)$ is not a central subloop of $Q$, then $\inn Q$ is not an abelian group.
If $A(Q) \le Z(Q)$, then $\inn Q$ is abelian if and only if
$$[z,u][z,v]^u = [z,v][z,u]^v \text{ or, equivalently, \,}
[z,vu][z,v,u] = [z,uv][z,u,v],$$
for all $u,v,z \in Q$.
\end{prop}
\dem If $\inn Q$ is abelian, then $A(Q) \le Z(Q)$, by \lref{72}. Assume
$A(Q) \le Z(Q)$. Then $Q/N$ and $A(Q)$ are boolean groups, 
and the associators are cyclically invariant, by \lref{42}.
In light of \lref{72} it is clear that if $A(Q) \le Z(Q)$,
then $\inn Q$ is abelian if and only if the mappings $z\mapsto z[z,u]$
and $z\mapsto z[z,v]$ commute, for all $u,v \in Q$. This gives us
the equality
$$z[z,u][z[z,u],v] = z[z,v][z[z,v],u], \text{ for all } u,v,z \in Q.$$
By \lref{41}, $[z,u][z[z,u],v] = [z,u][z,v]^{[z,u]}[[z,u],v] = [z,v][z,u]
[[z,u],v]$. Hence $$[z,v][z,u][[z,u],v] = [z,u][z,v][[z,v],u]$$ is a condition
that expresses the commutativity of the above mappings.

We have $[x,y] = [y,x]\m$, since $[x,y] \in  N(Q)$, for all $x,y \in Q$.
From \lref{41} we hence get the general equality
$$ [z,xy] =[z,y][z,x]^y[z,y,x], \text{ for all } x,y,z \in Q.$$

Now, $[z,v][z,u][[z,u],v] = [z,v][z,u][z,u]\m [z,u]^v = [z,v][z,u]^v
= [z,uv][z,v,u]$, and the rest is clear.
\edem

For a loop $Q$ one can define $Q'$ in a similar way as in groups, i.e.~
as the least normal subloop $S$ such that $Q/S$ is an abelian group.
If $N(Q) \unlhd Q$ and $Q/N(Q)$ is abelian, then clearly $Q' \le N(Q)$.
This is so in every Buchsteiner loop $Q$, and hence $Q'$ has to
be always a group in such loops. Note, that \lref{41}
can be used to see that $Q'$ coincides with the subloop generated
by all associators $[x,y,z]$ and commutators $[x,y]$.

\begin{lem}\label{74}
Let $Q$ be a Buchsteiner loop with $\inn Q$ abelian. Then $Q'$
is abelian as well.
\end {lem}
\dem Consider $x,y_1,y_2,z \in Q$ and put $y = y_1y_2$. Then $[z,y_1y_2]$
is equal to
$[z,y_2][z,y_1]^{y_2} [z,y_2,y_1]$, and so $[z,y][z,x]^y = [z,x][z,y]^x$
gives
$$[z,y_2][z,y_1]^{y_2}[z,x]^{y_1y_2} = [z,x][z,y_2]^x[z,y_1]^{y_2x}
= [z,y_2][z,x]^{y_2}[z,y_1]^{y_2x}.$$
We also have $[z,y_1]^{y_2}[z,x]^{y_1y_2} = [z,x]^{y_2}[z,y_1]^{xy_2}$,
and that means that $[z,y_1]^{y_2x} = [z,y_1]^{xy_2}$. In other words,
$[u_1,u_2]^{vw} = [u_1,u_2]^{wv}$ for all $u_1,u_2,v,w\in Q$. That is the
same as $[u_1,u_2]^{[v,w]} = [u_1,u_2]$, and so $[u_1,u_2][v,w]
= [v,w][u_1,u_2]$. Each commutator and each associator thus commutes
with every commutator and with every associator. The group $Q'$
is hence abelian.
\edem

Let $Q$ be a Buchsteiner loop. Then $Q/A(Q)$ acts upon $N(Q)$,
and so also upon $Q'$. If $Q'$ is abelian, then we get
an action of $Q/Q'$ upon $Q'$, and so we get an action
of an abelian group upon an abelian group. We shall
use this fact freely in the following lemma, understanding
that $T_x\m([u,v]) = [u,v]^x$ means in fact the action of $xQ'$
upon $[u,v]$, and so $[u,v]^{yx} = [u,v]^{yx}$, for all $x,y \in Q$.

\begin{lem}\label{75}
Let $Q$ be a Buchsteiner loop with $A(Q) \le Z(Q)$ and $Q'$ abelian
that is generated by a set $X$. If $[z,y][z,x]^y = [z,x][z,y]^x$ holds
for all $x,y,z \in X$, then it holds for all $x,y,z \in Q$.
\end{lem}
\dem Let us have $y = y_1y_2$, where $y_1,y_2 \in Q$. We can express $[z,y_1y_2]$
as $[z,y_2][z,y_1]^{y_2}[z,y_2,y_1]$, and we see that $[z,y][z,x]^y$ equals
$[z,x][z,y]^x$ if and only if $[z,y_2][z,y_1]^{y_2}[z,x]^{y_1y_2}$
equals $[z,x][z,y_2]^{x}[z,y_1]^{y_2x}$. The equality thus takes place
if and only if $a = b^{y_2}$, where $a=([z,x][z,y_2]^x)\m[z,y_2][z,x]^{y_2}$
and $b = [z,y_1]^x([z,x]^{y_1}[z,y_1])\m[z,x]$. Hence
\begin{align*}
a = b^{y_2} &\ \Leftrightarrow \ [z,y][z,x]^y = [z,x][z,y]^x, \\
a = 1 &\ \Leftrightarrow \ [z,y_2][z,x]^{y_2} = [z,x][z,y_2]^x, \text{ and}\\
b= 1 &\ \Leftrightarrow \ [z,y_1][z,x]^{y_1} = [z,x][z,y_1]^x.
\end{align*}
(Note that we have been using the commutativity of $Q'$ when expressing
the condition $b = 1$). If two of conditions $a = b^{y_2}$, $a = 1$ and
$b = 1$ are true, then the third one is true as well. From that we see
that if $[z,u_i][z,x]^{u_i} = [z,x][z,u_i]^x$ holds for $i \in \{1,2\}$,
then $[z,u][z,x]^u = [z,x][z,u]^x$ for every $u \in \{u_1u_2, u_1/u_2,
u_1\bs u_2\}$.

The case $z = z_1z_2$ is similar. We have $[z_1z_2,y][z_1z_2,x]^y
= [z_1z_2,x][z_1z_2,y]^x$ if and only if $[z_1,y]^{z_2}[z_2,y][z_1,x]^{z_2y}
[z_2,x]^y$ equals $[z_1,x]^{z_2}[z_2,x][z_1,y]^{z_2x}[z_2,y]^x$. The equality
takes place if and only if $a^{z_2} = b$, where $a = [z_1,y][z_1,x]^y
([z_1,x][z_1,y]^x)\m$ and $b = [z_2,x][z_2,y]^x([z_2,y][z_2,x]^y)\m$.
The rest is clear.
\edem

\section{Construction of a Buchsteiner loop of order 128} \label{8} 

The purpose of this section is to show that there exist proper
Buchsteiner loops with $\inn Q$ abelian. Such loops cannot be
of nilpotency class two, since then they would be conjugacy
closed, by \cref{44}.

We shall be constructing a loop $Q$ with $N(Q) = Q'$,
$Q/Q' \cong C_2 \times C_2 \times C_2$ and $Q' \cong C_4 \times
C_2 \times C_2$. The loop will be defined by a traditional
method upon the set $B\times N$, where $B \cong Q/Q'$ is written
multiplicatively and $N \cong N(Q)$ additively, $B$ acts
multiplicatively upon $N$, and
$$(u,a)\cdot (v,b) = (uv,\, \theta(u,v) + va + b) \text{ for all }
u,v \in B \text{ and } a,b \in N,$$
where the factor system (2-cocycle) $\theta: B\times B \to N$ is
defined in such a way that $\theta(u,1) = 0 = \theta(1,u)$ for
all $u \in B$.

For loops defined in this way (with both $B$ and $N$ abelian)
one can easily compute the associator $[(u,a),(v,b),(w,c)]$
as
$$(1,\, \theta(uv,w)+w\theta(u,v) - \theta(u,vw) - \theta(v,w)),$$
which means that $1 \times B$ is always contained in the
nucleus (and thus also in the commutant).

Assume that $B$ is generated by $e_i$, $1 \le i \le 3$,
and that $N$ is generated by $h$, an element of order four,
and by a subgroup $\{0,c_1,c_2,c_3\}$ that is isomorphic
to $C_2 \times C_2$. We shall write $2h$ sometimes as $d$,
and so $-h = d+h$.

Define a (multiplicative) action of $B$ upon $N$ by
$$ e_i h = h + d,\ e_i c_i = c_i \text{ and } e_i c_j = c_j + d,$$
for all $i,j \in \{1,2,3\}$, $i \ne j$. Clearly, $ud = d$,
for all $u \in B$, and $B$ acts trivially upon $N/D$,
$D = \{0, d\}$.

We shall define $\theta: B \times B \to N$ as a sum, with
$\theta(u,v) = \eta(u,v) + \delta(u,v)d$ for all $u,v \in B$, where
$\eta: B \times B \to N$, and $\delta: B \times B \to \{0,1\}$.
Now,
$$\eta(\prod e_i^{\alpha_i}, \prod e_i^{\beta_i}) =
\sum\alpha_i\beta_{i-1}h + \sum(\alpha_i\beta_i +
\alpha_{i-1}\beta_{i+1})c_i,$$
where $\alpha_i,\beta_i\in \{0,1\}$, and the indices are
computed modulo three. Furthermore,

$$\delta(\prod e_i^{\alpha_i}, \prod e_i^{\beta_i})
= \sum \alpha_i \alpha_{i+1}\beta_i + \sum (\alpha_i +
\alpha_{i-1})\beta_i\beta_{i+1},$$
again for all $\alpha_i,\beta_i\in \{0,1\}$ (the indices are computed
modulo three, and the expression is computed modulo two).

The mapping $\eta$ is defined so that $\eta(e_i,e_i) = c_i$,
$\eta(e_i,e_{i+1}) = 0$ and $\eta(e_i,e_{i-1})=h+c_{i+1}$.

\begin{lem}\label{81}
Assume  $\alpha_i,\beta_j\in \{0,1\}$, $i,j\in\{1,2,3\} 3$. Then
$$ \eta(\prod e_i^{\alpha_i}, \prod e_j^{\beta_j}) =
\sum \alpha_i\beta_j\eta(e_i,e_j).$$
\end{lem}
\dem We have $\sum \alpha_i\beta_j\eta(e_i,e_j) =
\sum \alpha_i\beta_ic_i + \sum \alpha_i\beta_{i-1}
(h+c_{i+1}) = \sum \alpha_i\beta_{i-1}h +
\sum(\alpha_i\beta_i + \alpha_{i-1}\beta_{i+1})c_i.$
\edem

\lref{81} seems to suggest that $h(uv,w) = h(u,w) + h(v,w)$
and $h(w,uv) = h(w,u) + h(w,v)$, for all $u,v,w \in B$.
However, none of these two equalities holds in general.
The reason is that $B$ is of exponent two and $\eta(e_i,
e_{i-1})$ is an element of order four. Nevertheless,
it is not difficult to compute the correction terms.
For that we shall use $\oplus$ as the addition modulo $2$
upon $\{0,1\}$. Note that for $\alpha,\beta \in \{0,1\}$
we always have $\alpha \oplus \beta = \alpha + \beta -
2\alpha\beta$, where the addition on the right hand side
is that of integers.

\begin{lem}\label{82}
Let $u = \prod e_i^{\alpha_i}$, $v = \prod e_i^{\beta_i}$
and $w = \prod e_i^{\gamma_i}$ be elements of $B$.
Then
\begin{align*}
\eta(uv,w) - \eta(u,w) - \eta(v,w) &= (\sum \alpha_{i+1}\beta_{i+1}\gamma_i)d,
\text{ and}\\
\eta(u,vw) - \eta(u,v) - \eta(v,w) &= (\sum \alpha_{i+1}\beta_i\gamma_i)d.
\end{align*}
\end{lem}

\dem Set first $\lambda_i = \alpha_i + \beta_i - 2\alpha_i\beta_i = \alpha_i \oplus\beta_i$, $1 \le i \le 3$. Then $uv = \prod e_i^{\lambda_i}$, and
$\eta(uv,w) = \sum \lambda_i \gamma_{i-1} h + \sum (\lambda_i\gamma_i
+ \lambda_{i-1}\gamma_{i+1})c_i = \sum(\alpha_i\beta_i\gamma_{i-1})d
+ \sum(\alpha_i\gamma_{i-1} + \beta_i\gamma_{i-1})h +
\sum(\alpha_i\gamma_i + \beta_i\gamma_i + \alpha_{i-1}\gamma_{i+1}
+ \beta_{i-1}\gamma_{i+1})c_i$, and that makes the former equality
clear. For the latter one proceed similarly, set $\nu_i
= \beta_i + \gamma_i - 2\beta_i\gamma_i$ and note that
$\alpha_i\nu_{i-1}h = (\alpha_i\beta_i + \alpha_i\gamma_i)h
+ \alpha_i\beta_{i-1}\gamma_{i-1}$.
\edem

To be able to utilize \lref{82} in the computation of the associator
we need to be able to express the difference of $w\eta(u,v)$ and $\eta(uv)$.
This is the content of the next lemma.

\begin{lem}\label{83}
Let $u = \prod e_i^{\alpha_i}$, $v = \prod e_i^{\beta_i}$
and $w = \prod e_i^{\gamma_i}$ be elements of $B$, and
let $x = \lambda h + \sum \rho_jc_j$ be an element of $N$.
Then $wx - x = \sum \gamma_i(\lambda + \rho_{i-1} + \rho_{i+1})d$
and
$$ w\eta(u,v) - \eta(u,v) = \sum (\alpha_{i-1}\beta_{i-1} +
\alpha_{i+1}\beta_{i+1} + \alpha_{i-1}\beta_{i+1})\gamma_id.$$
\end{lem}
\dem First note that the formula for $wx -x$ is defined correctly.
Indeed, set $\rho'_i = \rho_i + \rho_3$ for $i \in \{1,2\}$
and set $\rho'_3 = 0$. Then $(\rho_{3-1} + \rho_{3+1})d =
(\rho'_{3-1} + \rho'_{3+1})d$ and for $i \in \{1,2\}$ we
get $\rho_{i-1} + \rho_{i+1} = \rho'_{i-1} + \rho'_{i+1}$.
The mapping $x\mapsto \sum \gamma_i(\lambda + \rho_{i-1}+\rho_{i+1})d$
thus yields an endomorphism of the abelian group $N$.
The mapping $x \mapsto wx -x$ is also such an endomorphism, and
hence it suffices to verify that both endomorphisms agree for
$x = h$ and $x = c_j$. Howewer, that comes immediately from the
definition of the action of $B$ upon $N$.

We have to apply the endomorphism to $x = \eta(u,v)$, which
means that $\lambda = \sum\alpha_i\beta_{i-1}$ and $\rho_j = \alpha_j\beta_j +
\alpha_{j-1}\beta_{j+1}$. Each $\gamma_j d$ is hence multiplied
by $\alpha_j\beta_{j-1} + \alpha_{j-1}\beta_{j+1} + \alpha_{j+1}\beta_j
+ \alpha_{j-1}\beta_{j-1} + \alpha_{j+1}\beta_j + \alpha_{j+1}\beta_{j+1}
+ \alpha_j\beta_{j-1}$, and that is equal modulo 2 to $\alpha_{j-1}\beta_{j+1}
+ \alpha_{j-1}\beta_{j-1} + \alpha_{j+1}\beta_{j+1}$, for all $j \in
\{1,2,3\}$.
\edem

\begin{cor}\label{84}
Let $u = \prod e_i^{\alpha_i}$, $v = \prod e_i^{\beta_i}$
and $w = \prod e_i^{\gamma_i}$ be elements of $B$.
Then $\eta(uv,w) + w\eta(u,v) - \eta(u,vw) - \eta(v,w)$
is equal to $\sum(\alpha_{i-1}\beta_{i-1} + \alpha_{i-1}\beta_{i+1}
+ \alpha_{i+1}\beta_i)\gamma_id$.
\end{cor}

\begin{lem}\label{85}
Let $u = \prod e_i^{\alpha_i}$, $v = \prod e_i^{\beta_i}$
and $w = \prod e_i^{\gamma_i}$ be elements of $B$. Then
$\delta(u+v,w) + \delta(u,v) + \delta(u,v+w) + \delta(v,w)$
is modulo 2 equal to $\sum(\alpha_{i-1}\beta_{i-1} +
\alpha_{i-1}\beta_{i+1} + \alpha_{i+1}\beta_i + \alpha_{i+1}
\beta_{i-1})\gamma_i$.
\end{lem}
\dem By definition,
$$\delta(u+v,w) = \sum(\alpha_i + \beta_i)(\alpha_{i+1} + \beta_{i+1})\gamma_i
+ \sum(\alpha_i + \alpha_{i-1} + \beta_i + \beta_{i-1})\gamma_i\gamma_{i+1},$$
which is clearly equal to $\delta(u,w) + \delta(v,w) +
\sum(\alpha_i\beta_{i+1} + \alpha_{i+1}\beta_i)\gamma_i.$
Similarly,
$$ \delta(u,v+w)= \sum\alpha_i\alpha_{i+1}(\beta_i + \gamma_i)
+ \sum (\alpha_{i-1}+\alpha_i)(\beta_i + \gamma_i)(\beta_{i+1}+\gamma_{i+1})$$
is equal to $\delta(u,v) + \delta(u,w) + \sum(\alpha_{i-1} +
\alpha_i)(\beta_i\gamma_{i+1} + \beta_{i+1}\gamma_i)$, and the latter
sum can be clearly expressed also as $\sum (\alpha_{i-1}\beta_{i+1}
+ \alpha_i\beta_{i+1} + \alpha_{i+1}\beta_{i-1} +
\alpha_{i-1}\beta_{i-1})\gamma_i$. The rest is obvious.
\edem

\begin{prop}\label{86}
The loop $Q$ is a Buchsteiner loop that is not conjugacy closed.
It is of nilpotency class three and its inner mapping group is abelian.
The nucleus of $Q$ is equal to $1 \times N$ and coincides with
$Q'$, the centre is equal to $\{(1,0),(1,d)\}$ and coincides
with $A(Q)$. Finally, $Z(Q/Z(Q)) = N(Q)/Z(Q)$.
\end{prop}
\dem
Let $u = \prod e_i^{\alpha_i}$, $v = \prod e_i^{\beta_i}$
and $w = \prod e_i^{\gamma_i}$ be elements of $B$. The associator
is given by $\theta(uv,w)+w\theta(u,v) - \theta(u,vw) - \theta(v,w))$
which equals the sum $\eta(uv,w)+w\eta(u,v) - \eta(u,vw) - \eta(v,w))$
and of $(\delta(uv,w)+\delta(u,v) + \delta(u,vw) + \delta(v,w)))d$
since $w\theta(u,v) = w\eta(u,v) + \delta(u,v)d$. From \cref{84}
and \lref{85} we hence get
$$[u,v,w] = (1,\sum \alpha_{i+1}\beta_{i-1}\gamma_i)d)
= (1, \sum \alpha_{i-1}\beta_i\gamma_{i+1})d).$$
Loop $Q$ has to be a Buchsteiner loop since  $[u,v,w] = [v,w,u] = [w,v,u]$
is a central element of exponent $2$, for all $u,v,w \in Q$.
If $e_i$ is identified with $(e_i,0)$, and $(1,a)$ with $a$, for all
$a \in N$, then we get $[e_1,e_2,e_3] = d$, $[e_1,e_3,e_2] = 0$,
and the other associator values can be computed by cyclic shifts and
by linearity.

One needs to multiply $(v,b) \cdot (u,b) = (vu,\theta(v,u) + ub + a)$
by $(0,\theta(u,v)-\theta(v,u) + (v-1)a + (1-u)b)$ to get
$(u,a)\cdot(v,b) = (uv,\theta(u,v) + va + b)$. Thus $[e_i,e_j]
= \theta(e_i,e_j) - \theta(e_j,e_i)$, and we get
$$[e_i,e_{i+1}]=h+d+c_{i-1} \text{ and } [e_i,e_{i-1}] = h + c_{i+1}.$$
Furthermore, $[e_ie_{i+1},e_{i-1}e_i] = \eta(e_ie_{i+1},e_ie_{i-1}) - \eta
(e_{i}e_{i-1},e_ie_{i+1}) = (d + c_i + c_{i+1} + c_{i-1}) - h = h$.
It is hence clear that $Q'$ is equal to $N(Q) = 1\times N$. To see
that $Z(Q/Z(Q)) = N(Q)/Z(Q)$ it remains to verify that $(e_1e_2e_2,0)$
does not commute with all elements of $Q$ modulo $Z(Q)$. However,
$\theta(e_1,e_1e_2e_3) - \theta(e_1e_2e_2,e_1) = (h + c_3) - (h + c_2 +d)
= d + c_1 \notin Z(Q)$.

To finish the proof we need to show that $[e_i,e_k]+ e_i[e_j,e_k] =
[e_j,e_k]+ e_j[e_i,e_k]$ for all $i,j,k \in \{1,2,3\}$, by \lref
{84}. The case $i = j$ is trivial, and so we can assume $j = i+1$,
by the symmetry of $i$ and $j$. If $k =i$, then $e_i[e_{i+1},e_i]
= e_i(h+c_{i-1}) = h+c_{i-1} =[e_{i+1},e_i]$. If $k = j = i+1$,
then $[e_i,e_{i+1}] = h + d + c_{i-1} = e_{i+1}[e_i,e_{i+1}]$.
Finally, let us have $k = i-1$. Then $[e_i,e_{i-1}]+e_i[e_{i+1},e_{i-1}]
= (h+c_{i+1}) + e_i(h+d+c_i) = d + c_{i-1} = (h+d+c_i) + e_{i+1}(h+c_{i+1})
= [e_{i+1},e_{i-1}] + e_{i+1}[e_i,e_{i-1}]$.
\edem


\begin{thebibliography}{9}
\bibitem{con} 
R. H. Bruck:  Contributions to the theory of loops,
\emph{Trans. Amer. Math. Soc.}
\textbf{60} (1946), 245--354.
\bibitem{hhb} 
H. H. Buchsteiner: O nekotorom klasse binarnych lup, \emph{Mat. Issled.}
\textbf{39} (1976), 54--66.
\bibitem{csd} 
P. Cs\"org\H o and A. Dr\'apal: Left conjugacy closed loops of nilpotency class two,
\emph{Results Math.} \textbf{47} (2005), 242--265.
\bibitem{csa} 
P. Cs\"org\H o: On connected transversals to abelian subgroups and 
loop theoretical consequences, \emph{Arch. Math. (Basel)}
\textbf{86}  (2006), 499--516.
\bibitem{csc} 
P. Cs\"org\H o: Abelian inner mappings and nilpotency class greater than two,
\emph{Europ.~J.~Comb.} \textbf{28} (2007), 858--868.
\bibitem{ccc} 
P. Cs\"org\H o and A. Dr\'apal: Loops that are conjugacy
closed modulo the center (submitted).
\bibitem{bcc} 
P. Cs\"org\H o and A. Dr\'apal: Buchsteiner loops and conjugacy closedness
(submitted).
\bibitem{buch} 
P. Cs\"org\H o, A. Dr\'apal and M. Kinyon: Buchsteiner loops
(submitted).
\bibitem{nuc} 
A. Dr\'apal, P. Jedli\v cka: On loop identities that can be obtained by
a nuclear identification (submitted).
\bibitem{dia} 
M. K. Kinyon, K. Kunen, J. D. Phillips: Diassociativity in conjugacy closed loops,
\emph{Comm. Algebra} \textbf{32} (2004), 767--786.
\bibitem{pac} 
M. K. Kinyon, K. Kunen: Power-associative, conjugacy closed loops,
\emph{J. Algebra} \textbf{304} (2006), 269--294.
\end{thebibliography}
\end{document}